\documentclass[12pt]{amsart}

\usepackage{amsmath,amsfonts,amssymb,amsthm}
\usepackage{graphics,graphicx,color}

\newtheorem{theorem}{Theorem}[section]
\newtheorem{lemma}[theorem]{Lemma}
\newtheorem{proposition}[theorem]{Proposition}
\newtheorem{corollary}[theorem]{Corollary}

\def \Rm {\mathbb R}

\def \A {\mathcal{A}}

\begin{document}

\title{Uniqueness and non-uniqueness in inverse radiative transfer}


\author{Plamen Stefanov}
\thanks{First author partly supported by  NSF FRG Grant No.~0554065}
\address{\hskip-\parindent
Plamen Stefanov\\Department of Mathematics\\
Purdue University\\
150 N. University Street\\
West Lafayette, IN 47907-2067} \email{stefanov@math.purdue.edu}
\author{Alexandru Tamasan}
\address{\hskip-\parindent
Alexandru Tamasan\\
Department of Mathematics\\
University of Central Florida\\
4000 Central Florida Blvd.\\Orlando, FL, 32816, USA}
\email{tamasan@math.ucf.edu}

\begin{abstract}
We characterize the non-uniqueness in the inverse problem for
  the stationary transport model, in which the absorption $a$ and the scattering coefficient $k$
  of the media are to be recovered from the albedo operator. We show that
  ``gauge equivalent'' pairs $(a,k)$ yield the same albedo operator,
  and that we can recover uniquely the class of the gauge equivalent pairs. We apply this result
  to show unique determination of the media when the absorption $a$ depends on the
  line of travel through each point while scattering $k$ obeys a symmetry
  property.  Previously known results concerned directional independent absorption $a$.
\end{abstract}

\maketitle



\pagestyle{myheadings}
\markboth{P. Stefanov and A. Tamasan}{Uniqueness and non-uniqueness
in inverse radiative transfer}


\section{Introduction}
\label{sec:intro}
%
This paper considers the problem of recovering the absorption and
scattering properties of a bounded, convex medium
$\Omega\subset\Rm^n$, $n\geq 3$ from the spatial-angular
measurements of the density of particles at the boundary
$\partial\Omega$. Provided that the particles interact with the
medium but not with each other, the radiation transfer in the
steady-state can be modeled  by the transport equation
\begin{equation}\label{transport_eq}
-\theta\cdot\nabla
u(x,\theta)-a(x,\theta)u(x,\theta)+\int_{S^{n-1}}k(x,\theta',\theta)u(x,\theta')d\theta'=0,
\end{equation}
for $x\in\Omega$ and $\theta\in S^{n-1}$; see, e.g.
\cite{caseZweifel67, RS}. The function $u(x,\theta)$ represents the
density of particles at $x$ traveling in the direction $\theta$,
$a(x,\theta)$ is the absorption coefficient at $x$ for particles
traveling in the direction of $\theta$, and $k(x,\theta',\theta)$
is the scattering coefficient (or the collision kernel) which
accounts for particles from an arbitrary direction $\theta'$ which
scatter in the direction of travel $\theta$.

The medium is probed with the given radiation
\begin{equation}\label{bd_condition_in}
u|_{\Gamma_-}=f_-
\end{equation}
and the exiting radiation
\begin{equation}\label{bd_condition_out}
u|_{\Gamma_+}=:\A [f_-]
\end{equation}is detected; where
\begin{equation}
  \label{gammapm}
  \Gamma_\pm = \{(x,\theta)\in \partial\Omega\times S^{n-1}:\quad
    \pm \theta\cdot n(x)>0 \}
\end{equation}with $n(x)$ denoting the outer unit normal at a boundary point $x$.
The equation \eqref{bd_condition_out} defines the albedo operator
$\A$ which takes the incoming flux $f_-$ to the outgoing flux
$u|_{\Gamma_+}$ at the boundary.

In general, the boundary value problem \eqref{transport_eq} and
\eqref{bd_condition_in} may not be uniquely solvable but it has a
unique  solution for generic $(a,k)$; see \cite{stefanovUhlmann08}
and Proposition~\ref{pr_1} below. Unique solvability can also be
obtained under some physically relevant subcritical conditions like
\eqref{critic_CS} or \eqref{critic_DL} below. We assume that $(a,k)$
is such that the direct problem \eqref{transport_eq} and
\eqref{bd_condition_in} is well posed.

The inverse boundary value problem of radiative transfer is to
recover the absorption $a(x,\theta)$ and the scattering kernel
$k(x,\theta',\theta)$ from knowledge of the albedo operator $\A$.
One could also study the time-dependent version of
\eqref{transport_eq}, and then the kernel of $\A$ contains one
more variable that gives us more information. This problem has
been solved under some restrictive assumptions (e.g.  $k$ of a
special type or independent of a variable) in
\cite{anikonov75,anikonov84, anikonovBubnov88,bal00, larsen84,
Larsen, mcCormick92,  tamasan02, tamasan03}. In three or higher
dimensions, uniqueness and reconstruction results for general $k$
and $a=a(x)$ were established in \cite{choulliStefanov99}. The
approach there is based on the study of the singularities of the
fundamental solution of \eqref{transport_eq} (see also \cite{Bo}),
and the singularities of the Schwartz kernel of $\A$. Stability
estimates for $k$ of special type   were established in
\cite{R,jnWang99}; and recently, for general $k$, in
\cite{balJollivet08}.  Uniqueness and reconstruction results in a
Riemannian geometry setting, including recovery of a simple
metric, were established in \cite{macdowall04}. Similar results
for the time-dependent model were established in \cite{CS1}, and
in \cite{KLU} for the Riemannian case. In planar domains the work
in \cite{stefanovUhlmann03} shows stable determination of the
isotropic absorption and small scattering, and extensions to
simple Riemannian geometry is given in \cite{macdowall05}. Also in
two dimensional domains we point out that the recovery of $k$ is
only known under smallness conditions which are more restrictive
than what is needed to solve the direct problem; e.g. more
restrictive than \eqref{critic_CS} or \eqref{critic_DL} below. On
the other hand, in the time-dependent case, the extra variable
allows us to treat the planar case without such restrictions, see
\cite{CS1}. We also mention here the recent works
\cite{balLangmoreMonard08, L, LM}, in which the coefficients are
recovered  from angularly averaged measurements rather than from
knowledge of the whole the albedo operator $\A$.

The above mentioned results concern media with directionally
independent absorption property; for transport with variable
velocity $v$, which now belongs to an open subset of $\Rm^n$, the
absorption may depend on the speed $a=a(x,|v|)$ .

In general, in  media with an anisotropic absorption, the albedo
operator does not determine the coefficients uniquely. For
example, if $k\equiv 0$ the obstruction to unique determination
can be readily seen. The most one can recover from the albedo
operator $\A$ are the integrals $\int_{\Rm}a(x+t\theta,\theta)dt$,
see \cite{choulliStefanov99}. In \cite{balJollivet08}, they are
shown to be stably recovered independently of $k$. In other words,
for each fixed direction $\theta\in S^{n-1}$, we know the integral
of the map $a(\cdot,\theta)$ over the parallel lines in the
direction of $\theta$. This is insufficient data, since one can
smoothly change the $x$-variable in the direction of $\theta$
while preserving the integral; see also the theorem below.

The non-uniqueness described above and the need to assume that $a$
is isotropic left the theory of the inverse radiative problem in a
somewhat unsatisfactory state. It was not clear whether the
uniqueness failed if $k\not\equiv0$, neither it was known what
information about $(a,k)$ can still be recovered. The purpose of
this work is to fill this gap. We show that a certain family of
``gauge transformations'' of $(a, k)$ does not change the albedo
operator $\A$; and that, given $\A$, one can recover uniquely the
class of gauge equivalent pairs $(a,k)$. The recovery of the gauge
equivalent class is explicit, as the recovery of $a=a(x)$ and $k$ in \cite{choulliStefanov99} is explicit.

\section{Main Results}
The pair $(a,k)$ is
assumed to satisfy the {\em admissibility} condition
\begin{equation}\label{atten_bd}
\sup_{(x,\theta)\in\Omega\times S^{n-1}}\left(|a(x,\theta)|+
\int_{S^{n-1}}|k(x,\theta,\theta')|d\theta'\right)<\infty.
\end{equation}
Let ${T}$ be the operator defined by the l.h.s. of
\eqref{transport_eq}. For $(a,k)$ admissible, the second and the
third terms of $T$ are bounded operators on $L^1(\Omega\times
S^{n-1})$, while the first term is unbounded. We view ${T}$ as an
unbounded operator on $L^1(\Omega\times S^{n-1})$ with the domain
\[
D({T})= \{f\in L^1(\Omega\times S^{n-1});\; \theta\cdot\nabla f\in L^1(\Omega\times S^{n-1}),\; f|_{\Gamma_-}=0\},
\]
see also \cite{choulliStefanov99}. Since the direct problem
\eqref{transport_eq} and \eqref{bd_condition_in} can always be
reduced to a non-homogeneous problem with a homogeneous boundary
condition, invertibility of ${T}$ implies well-posedness. We say
that the direct problem is well posed, if ${T}^{-1}$ exists as a
bounded operator.

As an example, we have the following two {\em subcritical}
conditions that yield well-posedness; see, e.g.,
\cite{balJollivet08,choulliStefanov99, dlen6, mokhtar, RS} and
Proposition~\ref{pr_1} below. Either
\begin{equation}\label{critic_CS}
\sup_{(x,\theta)\in\Omega\times S^{n-1}}
\left|\tau(x,\theta)\int_{S^{n-1}}k(x,\theta,\theta') d\theta'
\right|<1,
\end{equation}
where $\tau(x,\theta)$ is the total free path of a particle at $(x,\theta)$, see the beginning of the next section, or
\begin{equation}\label{critic_DL}
a(x,\theta)-\int_{S^{n-1}}k(x,\theta,\theta') d\theta'\geq
0,~~a.e.~~\Omega\times S^{n-1}.
\end{equation}
 both $a$ and $k$ outside $\Omega$ to


We start with a simple observation, that seems to be new. We will
show that there is non-uniqueness even if $k\not\equiv 0$, too.
Let $\phi(x,\theta)>0$ be such that $\phi=1$ on $\Gamma :=
\partial\Omega\times S^{n-1}$. Set $\widetilde a=a-\theta\cdot\nabla_x
\log \phi$. Then we can rewrite \eqref{transport_eq} as
\[
-\theta\cdot\nabla_x u- \phi^{-1}(\theta\cdot\nabla_x \phi )u
-\widetilde a u+ \int_{S^{n-1}}
k(\cdot,\theta',\theta)u(\cdot,\theta')d\theta'=0.
\]
Multiply by $\phi$ to get
\[
-\theta\cdot\nabla_x (\phi u) -\widetilde a\phi u+ \phi\int_{S^{n-1}}
k(\cdot,\theta',\cdot)u(\cdot,\theta')d\theta'=0.
\]
The function $\widetilde u = \phi u$ thus solves
\eqref{transport_eq} with $(a,k)$ replaced by $(\widetilde a,
\widetilde k)$, where those two pairs are related by the ``gauge transformation''
\begin{equation}\label{s01}
\widetilde a=a-\theta\cdot\nabla_x \log \phi, \quad \widetilde
k(x,\theta',\theta) =
\frac{\phi(x,\theta)}{\phi(x,\theta')}k(x,\theta',\theta).
\end{equation}
Since $\phi=1$ on $\Gamma$, the boundary data do not change.
Therefore, if $\widetilde\A$ is the albedo operator corresponding
to the pair $(\widetilde a, \widetilde k)$, then $\mathcal{A} =
\widetilde\A$.

Our main result is that this is the only obstruction to non-uniqueness.

\begin{theorem}\label{non_uniqueness}
Let $(a,k)$ and $(\widetilde a,\widetilde k)$ be two admissible
pairs for which the direct problem is well-posed, and let $\A$ and
$\widetilde \A$ be the corresponding albedo operators.
Then $\A = \widetilde A$ if and only if there exists a positive
$\phi\in L^\infty(\Omega\times S^{n-1})$ with $\theta\cdot\nabla_x
\phi(x,\theta)\in L^\infty(\Omega\times S^{n-1})$ and $\phi= 1$ for
$x\in \partial\Omega$, so that \eqref{s01} hold.
\end{theorem}

The proof of the theorem is based on the analysis of the singularities of the Schwartz kernel of $\A$, as in \cite{choulliStefanov99}.

Theorem~\ref{non_uniqueness} allows us to obtain a few new
uniqueness results under additional conditions. One of them
concerns the case where  we have anisotropic media with absorption
$a(x,\theta)$ that depends on the line of travel through each
point (but not on the direction):
\begin{equation}\label{sym_atten}
a(x,\theta)=a(x,-\theta),~~x\in\Omega,~~\theta\in S^{n-1},
\end{equation}
and a scattering coefficient $k> 0$ satisfying the following symmetry condition
\begin{equation}\label{sym_scat}
k(x,\theta,\theta')=k(x,\theta',\theta),~~x\in\Omega,~~\theta,\theta'\in
S^{n-1}.
\end{equation}

\begin{corollary}\label{symmetric_media} Let $(a,k)$, $(\widetilde{a},\widetilde{k})$ be two admissible
and subcritical pairs which yield the same albedo operator. Assume
that $k$ and $\widetilde{k}>0$ satisfy \eqref{sym_scat}.

(i) Then $k=\widetilde{k}$ and $a=\tilde a+\theta\cdot\nabla v(x)$ for some function $v(x)$ vanishing on $\partial\Omega$. In particular, one can recover the total absorption at a.e.\ $x$, i.e.,  $\int a(x,\theta)d\theta= \int\widetilde a(x,\theta)d\theta$.

(ii) If, in addition,  $a$ and $\widetilde{a}$ satisfy
\eqref{sym_atten}, then $a=\widetilde{a}$.
\end{corollary}

Note that any two pairs as in (i) yield the same albedo operator, so this is the most we can say in this case.
The symmetry assumption \eqref{sym_scat} occurs
naturally in some models of Optical Tomography, where the scattering of light in a
tissue depends on the angle between the two directions:
$k(x,\theta,\theta')=k(x,\theta\cdot\theta').$

One can formulate and prove similar results in the case where the
velocity belongs to an open subspace of $\Rm^n$, i.e., the speed
can change, as in \cite{choulliStefanov99}. We restrict ourselves
to the fixed speed case ($|\theta|=1$) for the sake of simplicity
of the exposition. Also, the fixed speed model is the one that is
most often discussed in the literature.

The paper is organized as follows. Section \ref{preliminaries}
recalls some results from \cite{choulliStefanov99} that we use
later. In Section \ref{main} we prove Theorem \ref{non_uniqueness}
and its Corollary \ref{symmetric_media}. In section
\ref{reconstruction} we give the reconstruction formulae for
continuous $a$ and $k$ in the symmetric case
covered by Corollary \ref{symmetric_media}. Section \ref{remarks}
contains concluding remarks.

\section{Preliminaries}\label{preliminaries}


In this section we recall some results from
\cite{choulliStefanov99} recast to the one-speed velocities and
introduce notations.

Let $\tau_\pm(x,\theta)$ be the travel time it takes a particle at
$x$ to reach the boundary while moving in the direction of
$\pm\theta$ and define
$\tau(x,\theta)=\tau_-(x,\theta)+\tau_+(x,\theta)$. Since we work
with unit-speed velocities, note that
$\tau(x,\theta)\leq\mbox{diam}(\Omega)$. Let
$d\xi(x,\theta)=|n(x)\cdot\theta|d\mu(x)d\theta$
where $d\mu(x)$ is the induced Lebesgue measure on the boundary
and $d\theta$ is the normalized measure on the sphere. Also let
$\delta_{\{x\}}(x')$ represent the delta distribution with
respect to the boundary measure $d\mu(x')$ supported at
$x\in\partial\Omega$ and $\delta_{\{\theta\}}(\theta')$ represents
the delta distribution with respect to $d\theta$ centered at
$\theta\in S^{n-1}$.

Also, for $(x,\theta)\in\Omega\times S^{n-1}$, let $x_\theta^+$
denote the exiting point if traveling from $x$ in the
$\theta$-direction and $x_{\theta'}^-$ be the entrance point at
the boundary to reach the inside point $x$ by traveling in the
$\theta'$-direction, i.e.
\begin{equation}\label{x^+/-}x_\theta^+:=x+\tau_+(x,\theta)~~\mbox{and}~~x_{\theta'}^-:=x-\tau_-(x,\theta')\theta'.
\end{equation}

For the proposition below, we introduce the class of
\textit{regular} scattering kernels $k\in C(\bar\Omega, \;
L^\infty(S^{n-1}, L^1(S^{n-1})))$. Then the map $\int
k(x,\theta',\theta)\phi(\theta') d\theta'$ is  bounded on
$L^1(S^{n-1})$ continuously depending on $x$. Our notion of regular
$k$ is stronger than that in \cite{mokhtar}, and, in particular, it
allows us to use the results in there.

\begin{proposition}\label{pr_1}
The direct problem is well-posed, i.e., $T$ has a bounded inverse
on $L^1(\Omega\times S^{n-1})$, if either \eqref{critic_CS} or
\eqref{critic_DL} holds. Moreover, the direct problem is
well-posed for an open dense set of $(a,k) \in
L^\infty(\Omega\times S^{n-1}) \times C(\bar\Omega, \;
L^\infty(S^{n-1}, L^1(S^{n-1})))$.

\end{proposition}
\begin{proof}
We will first discuss the well-posedness under the conditions
\eqref{critic_CS} or \eqref{critic_DL}. Assume \eqref{critic_CS}
first.
By \cite[Proposition~2.3]{choulliStefanov99}, and since $\tau$ is
bounded, we get that $T^{-1}$ is bounded. The subcritical case
\eqref{critic_DL} is covered in \cite[Section 2]{balJollivet08}.

The generic statement is proven in $L^2(\Omega\times S^{n-1})$ for
$C^2$ coefficients in \cite{stefanovUhlmann08}. In the $L^1$ spaces
under consideration, we proceed in a similar way. Let $K$ be the
integral operator in \eqref{transport_eq} and $T_1=T-K$. Then
$KT_1^{-1}K$ is weakly compact in $L^1(\Omega\times S^{n-1})$, see  \cite{mokhtar}. Therefore, $(KT_1^{-1})^2$ is weakly
compact, and its square is compact. For a fixed $a$, consider the
family $\lambda k$, where $\lambda$ is a real parameter. By the
analytic Fredholm alternative in Banach spaces \cite{RV},
$\lambda\mapsto (I-(\lambda KT_1^{-1})^4)^{-1}$ is a meromorphic
family. This implies that $\lambda\mapsto (I-\lambda
KT_1^{-1})^{-1}$ is also meromorphic, and thus  $T^{-1}= (I-\lambda
KT_1^{-1})^{-1}T_1^{-1}$ exists for all but a discrete set of
$\lambda$'s. This shows that there is a dense set of pairs yielding
a well-posed problem. The fact that this set is also open follows
from a perturbation argument around each $(a,k)$, for which $T^{-1}$
is bounded; thus $(I-\lambda KT_1^{-1})^{-1}$ corresponding
to nearby pairs exists.
\end{proof}
Note that one can set $a(x,\theta)=a_0(x,\theta)+\int
k(x,\theta,\theta')d\theta'$, where  the integral represents the
attenuation due to the change of direction, while $a_0$ is the
absorption. Then one can prove in the same way that the direct
problem is well posed for generic $(a_0,k)$, and moreover, for any
$a_0$, this is true for generic $k$'s.

\begin{proposition}
Assume that the direct problem is well-posed. Then the
albedo operator $\A:L^1(\Gamma_-,d\xi)\to L^1(\Gamma_+,d\xi)$ is
bounded and its Schwartz kernel $\alpha$ is given by
$\alpha=\alpha_1+\alpha_2+\alpha_3$, where
\begin{align}
\alpha_1(x,\theta,x',\theta')=&e^{\int_0^{\tau_-(x,\theta)}a_\theta(x-t\theta)dt}
\delta_{\{x_\theta^-\}}(x')\delta_{\{\theta\}}(\theta')\label{alpha_1}\\
\alpha_2(x,\theta,x',\theta')=&\int_0^{\tau_-(x,\theta)}
e^{-\int_0^s a(x-t\theta,\theta)dt}e^{-\int_0^{\tau_-(x-s\theta,\theta')}a(x-s\theta-t\theta',\theta')dt}\label{alpha_2}\\
&\times
k(x-s\theta,\theta',\theta)\delta_{\{x-s\theta-\tau_-(x-s\theta,\theta')\theta'\}}(x')ds\nonumber\\
|n(x')\cdot \theta'|^{-1}\alpha_3\in&
L^\infty(\Gamma_-;L^1(\Gamma_+,d\xi)).\label{alpha_3}
\end{align}
\end{proposition}

This proposition is formulated in \cite{choulliStefanov99} under
the assumption that the system is subcritical, i.e., either
\eqref{critic_CS} or \eqref{critic_DL} holds. We remark that those
conditions are only used in the proof in the analysis of
$\alpha_3$, to guarantee that $T^{-1}$ exists in $L^1$; something
that we assume here.

Let $\phi\in C^\infty_0(B(0;1))$ with $0\leq\phi\leq 1$ and
$\phi\equiv 1$ near the origin be a cut-off function. Given
$\epsilon>0$ we define for $x,x'\in\Rm^n$ and $\theta,\theta'\in
S^{n-1}$
\begin{equation}\label{phi_epsilon}
\phi_\epsilon (x,\theta,x',\theta')= \phi\left(\frac{x_\theta^{-}
- x'}{\epsilon}\right)
\phi\left(\frac{\theta-\theta'}{\epsilon}\right).
\end{equation}
\begin{proposition}\label{ballistic}Assume that the direct problem is well posed.
Then the limit below holds in $L^1(\Gamma_+,d\xi)$
\begin{equation}\label{limit_1}
\lim_{\epsilon\to
0}\int_{\Gamma_-}\alpha(x,\theta,x',\theta')\phi_\epsilon(x,\theta,x',\theta')d\mu(x')d\theta
=e^{-\int_{-\infty}^{\infty}a(x+t\theta,\theta)dt}>0~~a.e.~~
\Gamma_+;\end{equation}
\end{proposition}

For linearly independent $\theta,\theta'\in S^{n-1}$, let
$\pi_{\theta,\theta'}(x)$ denote the projection  of $x$ onto the
plane through the origin spanned by $\theta$ and $\theta'$. Let
$\theta'_\perp$ be a unit vector in $\mbox{span}\{\theta,\theta'\}$
orthogonal to $\theta'$: $\theta'_\perp\cdot\theta'=0$. Let
$\varphi\in C^\infty_0(-1,1)$ with $0\leq\varphi\leq 1$ and
$\int_\Rm\varphi(t)dt=1$ be a cut-off function in $\Rm$ and $\phi$
be the cut-off function in $\Rm^n$ introduced earlier. Define the
test function
\begin{align}\label{test_funct_2}
\phi_{\epsilon,\delta}(x;\theta,\theta')=\frac{1}{\epsilon}
\varphi\left(\frac{x\cdot\theta'_\perp}{\epsilon\theta\cdot\theta'_\perp}\right)
\phi\left(\frac{x-\pi_{\theta,\theta'}(x)}{\delta}\right).
\end{align}
With $x_\theta^+$ and $x_{\theta'}^-$ given by \eqref{x^+/-} and
$\phi_{\epsilon,\delta}$ given by \eqref{test_funct_2}, we define in
$\Omega\times S^{n-1}\times S^{n-1}$
\begin{align}\label{test_2}
I_{\epsilon,\delta}(x,\theta',\theta):= \int_{\partial\Omega}&
\alpha(x_\theta^+,\theta,x',\theta')
\phi_{\epsilon,\delta}(x'-x_{\theta'}^-,\theta,\theta')d\mu(x').
\end{align}

\begin{proposition}\label{scatt_once}Assume that the direct problem is well posed. Then
the following limit holds in
$L^1(\Omega;L^1_{loc}(I^\nshortparallel))$
\begin{align}\label{limit_2}
I_{0,0}(x,\theta',\theta):=\lim_{\epsilon\to 0}\lim_{\delta\to
0}I_{\epsilon,\delta}(x,\theta',\theta)=e^{-\int_{-\infty}^0
a(x+t\theta',\theta')dt-\int_0^{\infty}a(x+t\theta,\theta)dt}k(x,\theta',\theta),
\end{align}
where $I^\nshortparallel:=\{(\theta,\theta')\in S^{n-1}\times
S^{n-1}: ~~\theta'\neq\pm\theta\}$.
\end{proposition}

\section{Proof of Theorem~\ref{non_uniqueness} and Corollary \ref{symmetric_media}}\label{main}

\begin{proof}[Proof of Theorem~\ref{non_uniqueness}]
From  Proposition \ref{ballistic}, we can recover
\begin{equation}\label{s02}
\int_{-\infty}^\infty a(x+s\theta,\theta)ds, \quad
(x,\theta)\in\Omega\times S^{n-1}.
\end{equation}
In particular, the integrals above for  $f = a-\widetilde a$
vanish:
\begin{align}\label{vanish}
\int_{-\infty}^\infty f(x+s\theta,\theta)ds=0, \quad\forall
(x,\theta)\in\Omega\times S^{n-1}.
\end{align}
The kernel of the linear transform \eqref{s02} is easy to
describe. Set
\begin{align}\label{beam}
v(x,\theta) = \int_{-\infty}^0 f(x+s\theta,\theta)ds, \quad
(x,\theta)\in\Omega\times S^{n-1}.
\end{align}
Then $v\in L^\infty(\Omega\times S^{n-1})$ with
$\theta\cdot\nabla_x v=f\in L^\infty(\Omega\times S^{n-1})$ and
$v=0$ on $\Gamma_-$.  From \eqref{vanish} we get also that $v=0$ on
$\Gamma_+$, therefore $v(x,\theta)=0$ for $x\in\partial\Omega$ and
a.e. $\theta\in S^{n-1}$.

Set
\begin{align}\label{gauge}
\phi(x,\theta) = e^{v(x,\theta)}
\end{align}
 to get $a-\widetilde a= \theta\cdot\nabla_x \log \phi$ as
claimed. This shows the first part of \eqref{s01}.

Now, from Proposition \ref{scatt_once} we get
\begin{align*}
e^{ -\int_{-\infty}^0  a(x+t\theta',\theta')-  \int_0^{\infty}
a(x+t\theta,\theta)}   k(x,\theta',\theta)= e^{-\int_{-\infty}^0
\widetilde a(x+t\theta',\theta')-  \int_0^{\infty}  \widetilde
a(x+t\theta,\theta)} \widetilde   k(x,\theta',\theta).
\end{align*}
Using the first equality in \eqref{s01} that we already proved,
after a simple calculation, we derive
\begin{equation}\label{scatt_gauge}
\widetilde k(x,\theta',\theta)
=\frac{\phi(x,\theta)}{\phi(x,\theta')}k(x,\theta',\theta).
\end{equation}
\end{proof}

\begin{proof}[Proof of Corollary~\ref{symmetric_media}]
Next we use the characterization above to prove uniqueness in the
symmetric case. By swapping $\theta$ and $\theta'$ in
\eqref{scatt_gauge} we get
\begin{align*}
\widetilde k(x,\theta,\theta')
&=\frac{\phi(x,\theta')}{\phi(x,\theta)}k(x,\theta,\theta')=\frac{\phi(x,\theta')}{\phi(x,\theta)}k(x,\theta',\theta)\\
&=\frac{\phi(x,\theta')^2}{\phi(x,\theta)^2}\widetilde
k(x,\theta',\theta)=\frac{\phi(x,\theta')^2}{\phi(x,\theta)^2}\widetilde
k(x,\theta,\theta').
\end{align*}
For the second and fourth equality we used \eqref{sym_scat}, while
for the third one we used \eqref{scatt_gauge}. Since $\widetilde k$
does not vanish, we conclude that $\phi(x,\theta')=\phi(x,\theta)$
for all $\theta,\theta'\in S^{n-1}$. Therefore $\phi=\phi(x)$ is
independent of $\theta$ and applying \eqref{scatt_gauge} again we
conclude that $k=\widetilde k$.

So far we showed that $a(x,\theta)-\widetilde a(x,\theta)=\theta\cdot\nabla_x \phi(x)$.
If the symmetry relation \eqref{sym_atten} holds then
$a(x,\theta)-\widetilde a(x,\theta)=-\theta\cdot\nabla_x
\phi(x)$. Therefore $a-\widetilde a=0$.
\end{proof}

\section{Reconstruction formulas in the symmetric case}\label{reconstruction}
We showed in the previous section that under the symmetry
hypotheses \eqref{sym_scat} and \eqref{sym_atten} there is
uniqueness for $a$, $k$. The proof is constructive but it still
leads to the problem of recovering $a(x,\theta)$ first, up to
$\theta\cdot\nabla \log\phi$, from its under-determined X-ray
transform, see the l.h.s.\ of \eqref{vanish}. In this section, we
show that under the same assumptions, there is an explicit
reconstruction of a different, simpler type. Namely, we recover
the integrals in \eqref{vanish} first, that measure the total
absorption along straight lines,  but we are not trying to
determine $a$ from them. Then by \eqref{alpha_2} we recover
\[
e^{-\int_{-\infty}^0
a(x+t\theta',\theta')dt-\int_0^{\infty}a(x+t\theta,\theta)dt}k(x,\theta',\theta),
\]
see \eqref{broken_geodesic_1} below. We do not know the
attenuation term above because we know the attenuation along
straight but not broken lines. We can however swap $\theta$ and
$\theta'$ and multiply the results. Then we get $k^2$ multiplied
by the attenuation along the two lines through $x$ parallel to
$\theta$, and $\theta'$, respectively, and we know that
attenuation. This recovers $k$, and then we recover $a$.

We
need to strengthen the regularity assumptions on the coefficients to
\begin{align}\label{regularity+}
a\in C(\Omega\times S^{n-1}),~~ k\in C(\Omega\times S^{n-1}\times
S^{n-1}).
\end{align}
Firstly, we extend the limit in Proposition \ref{ballistic} valid
for maps in $\Gamma_+$ to maps in $\Omega\times S^{n-1}$. For
$(x,\theta)\in\Omega\times S^{n-1}$ and $\epsilon>0$, let us
denote
\begin{align}\label{test_1}
J_\epsilon(x,\theta):=\int_{\Gamma_-}\alpha(x_\theta^+,\theta,x',\theta')\phi_\epsilon
(x,\theta,x',\theta')d\mu(x')d\theta',
\end{align}where $x_\theta^+$ is given in \eqref{x^+/-}. We remark here that $J_\epsilon$
is constant with $x$ varying in the
direction of $\theta$, i.e.
$J_\epsilon(x+t\theta,\theta)=J_\epsilon(x,\theta)$ for all $t$.

\begin{lemma}\label{balistic} The limit below holds in $L^1(\Omega\times
S^{n-1};\tau(x,\theta)^{-1}dxd\theta)$
\begin{equation}\label{limit_11}
J_0(x,\theta):=\lim_{\epsilon\to
0}J_\epsilon(x,\theta)=\exp\left(-\int_{-\infty}^{\infty}a(x+t\theta,\theta)dt\right)>0~~a.e.~~
\Omega\times S^{n-1}.\end{equation}
\end{lemma}
\begin{proof}
We use the change of variable formula
$$\int_{\Omega\times S^{n-1}}f(x,\theta)dxd\theta=\int_{\Gamma_+}\int_0^{\tau_-(x^+,\theta)}f(x^+
-s\theta)dsd\xi(x^+,\theta)$$for
$f(x,\theta)=\left|J_{\epsilon}(x,\theta)-\exp\left(\int_{-\infty}^\infty
a(x+t\theta,\theta)dt\right)\right|\tau(x,\theta)^{-1}$ to
estimate
\begin{align*}
&\int_{\Omega\times
S^{n-1}}\left|J_{\epsilon}(x,\theta)-\exp\left(\int_{-\infty}^\infty
a(x+t\theta,\theta)dt\right)\right|\tau(x,\theta)^{-1}dxd\theta\\
&=\int_{\Gamma_+}\int_0^{\tau_-(x^+,\theta)}
\frac{\left|J_{\epsilon}(x^+-s\theta,\theta)-e^{\left(\int_{-\infty}^\infty
a(x^+-s\theta+t\theta,\theta)dt\right)}\right|}{\tau(x^+-s\theta,\theta)}{ds}d\xi(x^+,\theta)\\
&=\int_{\Gamma_+}
\left|J_{\epsilon}(x^+,\theta)-\exp\left(\int_{-\infty}^\infty
a(x^+ +t\theta,\theta)dt\right)\right|d\xi(x^+,\theta).
\end{align*}The last identity holds due to the fact that
$\tau(x^+-s\theta,\theta)=\tau(x^+,\theta)=\tau_-(x^+,\theta)$ and
that the integrand is constant in $s$. An application of
Proposition \ref{ballistic} finishes the proof.
\end{proof}

Following from Proposition \ref{scatt_once} for a.e.
$(x,\theta,\theta')\in\Omega\times I^\nshortparallel$ we get that
\begin{align}
e^{-\int_{-\infty}^0\label{broken_geodesic_1}
a(x+t\theta',\theta')dt-\int_0^{\infty}a(x+t\theta,\theta)dt}k(x,\theta',\theta)=I_{0,0}(x,\theta',\theta),\\
e^{-\int_{-\infty}^0
a(x+t\theta,\theta)dt-\int_0^{\infty}a(x+t\theta',\theta')dt}k(x,\theta,\theta')
=I_{0,0}(x,\theta,\theta').\label{broken_geodesic}
\end{align}

Using the symmetry assumption \eqref{sym_scat} on $k$, by
multiplication of \eqref{broken_geodesic_1} and
\eqref{broken_geodesic} we obtain
\begin{align}\label{k_intermediate}
e^{-\int_{-\infty}^\infty
a(x+t\theta',\theta')dt-\int_{-\infty}^{\infty}a(x+t\theta,\theta)dt}k^2(x,\theta',\theta)=
I_{0,0}(x,\theta',\theta)I_{0,0}(x,\theta,\theta').
\end{align}
Now using Lemma \ref{balistic} we recover
\begin{align}\label{k}
k(x,\theta',\theta)=
\left(\frac{I_{0,0}(x,\theta',\theta)I_{0,0}(x,\theta,\theta')}{J_0(x,\theta)J_0(x,\theta')}\right)^{\frac12}.
\end{align}

Using the formula \eqref{k} in \eqref{broken_geodesic} we get for
a.e. $(x,\theta,\theta')\in\Omega\times I^\nshortparallel$ that
\begin{align}
\int_{-\infty}^0a(x+t\theta,\theta)dt+\int_0^\infty
a(x+t\theta',\theta')dt=\frac12\log\left(\frac{I_{0,0}(x,\theta',\theta)}{I_{0,0}(x,\theta,\theta')J_0(x,\theta)J_0(x,\theta')}\right).
\end{align}
The continuity assumption \eqref{regularity+} imply that the
identity above extends pointwise in $\Omega\times S^{n-1}\times
S^{n-1}$. By looking at backscattering $\theta'=-\theta$ and
making one change of variables $t\leftrightarrow -t$ we obtain
\begin{align}\label{beams}
\int_{-\infty}^0a(x+t\theta,\theta)dt+\int_{-\infty}^0
a(x+t\theta,{-\theta})dt=\frac12\log\left(\frac{I_{0,0}(x,-\theta,\theta)}{I_{0,0}(x,\theta,-\theta)J_0(x,\theta)J_0(x,-\theta)}\right).
\end{align}
The left hand side is differentiable in $x$ in the direction of
$\theta$, hence so is the right hand side. By differentiating
\eqref{beams} in the direction of $\theta$ we obtain
\begin{align}
a(x,\theta)+a(x,{-\theta})=\frac12 \theta\cdot\nabla_x\log
\left(\frac{I_{0,0}(x,-\theta,\theta)}{I_{0,0}(x,\theta,-\theta)J_0(x,\theta)J_0(x,-\theta)}\right).
\end{align}
Since the absorption depends on the line through $x$, we have
$a(x,\theta)=a(x,{-\theta})$ and thus it can be recovered from the
formula:
\begin{align}\label{atten_reconst}
a(x,\theta)=\frac{1}{4}\theta\cdot\nabla_x\log
\left(\frac{I_{0,0}(x,-\theta,\theta)}{I_{0,0}(x,\theta,-\theta)J_0(x,\theta)J_0(x,-\theta)}\right).
\end{align}

\section{Remarks}\label{remarks}
\subsection{The isotropic absorption case}
The previously known uniqueness result \cite{choulliStefanov99}  for isotropic absorption follow from Theorem \ref{non_uniqueness} combined with the injectivity
of the X-ray transform. If $f=a-\widetilde{a}$ depends on the
position only, then \eqref{vanish} implies $f=0$. From the
definition \eqref{gauge} we get $\phi\equiv 1$, which by \eqref{s01} yields $k=\widetilde k$.

%
%
%
\subsection{Other conditions that imply uniqueness}
Assuming \eqref{sym_scat}, we can only recover $a$ up to
$\theta\cdot\nabla v(x)$ as above. For unique recovery, the
condition \eqref{sym_atten} suffices but it can be replaced by
something weaker. For example, we may require that $a$ is orthogonal
to all such functions w.r.t.\ some measure $d\mu(\theta)$. This is equivalent to
\[
\int_{S^{n-1}} \theta\cdot\nabla_x a(x,\theta)d\mu(\theta)=0.
\]
Note that the symmetry assumption \eqref{sym_atten}  implies
such a condition if $d\mu(\theta)$ is even.

{\bf Acknowledgement:} This work started during the Summer School on
Inverse Problems in Radiative Transfer at the University of
California at Merced. The authors thank Arnold Kim for his kind
hospitality.


\begin{thebibliography}{99}
\bibitem{anikonov75}D. S. Anikonov, {\em The unique determination of the coefficients
and the right side in a transport equation}, Differential
Equations {\bf 11}(1975), pp. 6--12.

\bibitem{anikonov84}\bysame {\em Multidimensional inverse problem for the transport
equation}, Differential Equations {\bf 20}(1984), pp. 608--614.

\bibitem{anikonovBubnov88}Yu. Anikonov and B. A. Bubnov, {\em Inverse problem of transport
theory}, Soviet Math. Dokl. {\bf 37}(1988), pp. 497--499.

\bibitem{Babovsky91}
{\sc H.~Babovsky}, {\em {Identification of Scattering Media from Reflected
  Flows}}, SIAM J. Appl. Math. {\bf 51}(1991), pp.~1674--1704.


\bibitem{bal00} {\sc G. Bal}, {\em Inverse problems for homogeneous transport
equations. II. The multidimensional case.}, Inverse Problems {\bf
16}(2000), pp.~ 1013--1028.

\bibitem{balJollivet08}
{\sc G.~Bal and A. Jollivet}, {\em {Stability estimates in
stationary inverse transport}}, preprint (2008).

\bibitem{balLangmoreMonard08}
{\sc G.~Bal, I. Langmore and F. Monard}, {\em {Inverse transport
with isotropic sources and angularly averaged measurements}},
Inverse Probl. Imaging {\bf 2}(2008), pp. 23-42.


\bibitem{Bo} {\sc A. Bondarenko}, {\em The structure of the fundamental
solution of the time-independent transport equation}, J. Math. Anal. Appl. {\bf 221}(1998), no. 2, 430--451.


\bibitem{caseZweifel67}
{\sc K.~M. Case and P.~F. Zweifel}, {\em {Linear Transport Theory}},
  Addison-Wesley series in nuclear engineering, Addison-Wesley, Reading, Mass.,
  1967.

\bibitem{CS1} {\sc M. Choulli and P. Stefanov}, {\em Inverse scattering
and inverse boundary value problems for the linear Boltzmann equation},
Comm. P.D.E. {\bf 21}(1996), 763--785.


\bibitem{choulliStefanov99}
{\sc M.~Choulli and P.~Stefanov}, {\em {An inverse boundary value
problem for the stationary transport}}, Osaka J. Math. {\bf
36}(1999), pp. 87--104.

\bibitem{dlen6}
{\sc R.~Dautray and J.-L. Lions}, {\em {Mathematical Analysis and
Numerical Methods for Science and Technology. Vol.6}}, Springer
Verlag, Berlin, 1993.

\bibitem{KLU}{\sc Y. Kurylev, M. Lassas, G. Uhlmann}, {\em  Rigidity of broken geodesic flow and inverse problems}, to appear in American Journal of Mathematics.

\bibitem{L}{\sc I. Langmore}, {\em The stationary transport problem with angularly averaged measurements}, Inverse Problems \textbf{24}(2008), no. 1, 015024.

\bibitem{LM} {\sc I. Langmore and S. McDowall}, {\ em Optical tomography
for variable refractive index with angularly averaged measurements}, reprint.

\bibitem{larsen84}
{\sc E. W. Larsen}, {\em Solution of multidimensional inverse
transport problems}, J. Math. Phys. {\bf 25}(1984), 131--135.

\bibitem{Larsen}{\sc E. W. Larsen}, {\em Solution of the three dimensional inverse transport problems},
Transport Theory and Stat. Phys. \textbf{17}(2\&3)(1988), 147--167.

\bibitem{mcCormick92}
{\sc N.~J. McCormick}, {\em {Inverse radiative transfer problems: a
review}}, Nucl. Sci. Eng. {\bf 112}(1992), pp.~185--198.

\bibitem{macdowall04}{\sc S. McDowall}, {\em Inverse problem for the transport equation in the presence of a Riemannian metric}, Pac. J. Math., 216 (2004), no.1, 107--129.


\bibitem{macdowall05}{\sc S. McDowall}, {\em Optical Tomography on Simple Riemannian
Surfaces}, Comm. PDE {\bf 30}, pp. 1379 -- 1400.

\bibitem{mokhtar}
{\sc M.~Mokhtar-Kharroubi}, {\em Mathematical Topics in Neutron
Transport Theory}, {World Scientific}, Singapore, 1997.

\bibitem{RS} {\sc M. Reed and B. Simon}, {\em Methods of Modern
Mathematical Physics}, Vol. 3, Academic Press, New York, 1979.

\bibitem{RV}{\sc M. Rivaric and I. Vidav}, {\em Analytic properties of the inverse $A(z)^{-1}$ of an analytic linear operator valued function $A(z)$}, Arch. Rat. Mech. Anal. \textbf{32}(4)(1969), 298--310.

\bibitem{R}{\sc V. Romanov}, {\em Stability estimates in problems of recovering the attenuation coefficient and the scattering indicatrix for the transport equation}, J. Inverse Ill-Posed Probl. \textbf{4}(1966)(4), 297--305.

\bibitem{stefanovUhlmann03}
{\sc P.~Stefanov and G.~Uhlmann}, {\em Optical tomography in two
dimensions}, Methods Appl. Anal. {\bf 10}(2003), pp.~1--9.

\bibitem{stefanovUhlmann08}
{\sc P.~Stefanov and G.~Uhlmann}, {\em An inverse source problem in optical molecular imaging}, to appear in Analysis and PDE.

\bibitem{tamasan02}
{\sc A.~Tamasan}, {\em An inverse boundary value problem in two-dimensional
  transport}, Inverse Problems {\bf 18}(2002), pp.~209--219.

\bibitem{tamasan03}\bysame, {\em Optical tomography in weakly anisotropic
scattering media}, Contemporary Mathematics {\bf 333}(2003), pp.
199--207.

\bibitem{jnWang99} J.-N. Wang, {\em Stability estimates of an inverse problem
for the stationary transport}, Ann. Inst. Henri Poincar\'e {\bf
70}(1999), pp.~473--495.

\end{thebibliography}
\end{document}